\documentclass[11pt,reqno]{amsproc}

\title[]{Estimates near the boundary for critical SQG}

\author{Peter Constantin}
\address{Department of Mathematics, Princeton University, Princeton, NJ 08544}
\email{const@math.princeton.edu}

\author{Mihaela Ignatova}
\address{Department of Mathematics, Temple University, Philadelphia, PA 19122}
\email{ignatova@temple.edu}

\usepackage[margin=1in]{geometry}
\usepackage{amsmath, amsthm, amssymb}
\usepackage{times}
\usepackage{color}
\usepackage{hyperref}

\newcommand{\pa}{\partial}
\newcommand{\la}{\label}
\newcommand{\fr}{\frac}
\newcommand{\na}{\nabla}
\newcommand{\be}{\begin{equation}}
\newcommand{\ee}{\end{equation}}
\newcommand{\ba}{\begin{array}{l}}
\newcommand{\ea}{\end{array}}
\newcommand{\Rr}{{\mathbb R}}

\newcommand{\beg}{\begin}

\newcommand{\D}{\Delta}
\renewcommand{\l}{\Lambda_D}
\date{today}
\begin{document}
\begin{abstract}
We obtain estimates near the boundary for the critical dissipative SQG equation in bounded domains, with the square root of the Dirichlet Laplacian dissipation.  We prove that global regularity up to the boundary holds if and only if a certain quantitative vanishing of the scalar at the boundary is maintained.   
\end{abstract}


\noindent

\maketitle

\section{Introduction}
The Surface Quasigeostrophic equation (SQG) of geophysical origin (\cite{held}) was proposed as a two dimensional model for the study of inviscid incompressible formation of singularities (\cite{c}, \cite{cmt}). The equation has been studied extensively. Blow up from smooth initial data is still an open problem, although  the original blow-up scenario of \cite{cmt} has been ruled out analytically (\cite{cord}) and numerically (\cite{cnum}. The addition of fractional dissipation  produces globally regular solutions if the power of the Laplacian is larger or equal than one half. When the linear dissipative operator is precisely the square root of the Laplacian, the equation is commonly referred to as the ``critical dissipative SQG'', or ``critical SQG''. The global regularity of solutions for critical SQG in the whole space or on the torus was obtained independently in \cite{caf} and \cite{knv} by very different methods.  Several subsequent proofs were obtained (see \cite{cvt} and references therein).

The critical SQG equation in bounded domains is given by

\be
\pa_t \theta + u\cdot\na\theta + \l \theta = 0
\la{sqg}
\ee
with
\be
u = \na^{\perp}\l^{-1}\theta.
\la{u}
\ee
Here $\Omega\subset \Rr^d$ is a bounded open set with smooth boundary, $\l$ is the square root of the Laplacian with vanishing Dirichlet boundary conditions, and $\na^{\perp} = J\na$ with $J$ an invertible antisymmetric matrix.  
The local existence and uniqueness of solutions of (\ref{sqg}) given in \cite{sqgb}
is
\beg{prop} \la{locsqg} Let $d=2$, and let $\theta_0\in H_0^1(\Omega)\cap H^{2}(\Omega) = {\mathcal{D}}(\l^2)$. There exists $T>0$ and  a unique solution of (\ref{sqg}) with initial datum $\theta_0$ satisfying
\be
\theta\in L^{\infty}(0,T; H_0^1(\Omega)\cap H^{2}(\Omega))\cap L^2\left(0,T; {\mathcal{D}}\left(\l^{2.5}\right) \right).
\la{loc}
\ee
\end{prop}
Local existence of solutions of the same type holds also for supercritical SQG in bounded domains (\cite{cn1}). A natural notion of strong solution is given by the condition
\be
\theta \in L^{\infty}(0,T; {\mathcal D}(\sqrt{\l}))\cap L^2(0,T; {\mathcal D}(\l))\la{strongsqg}.
\ee
The solutions of Proposition \ref{locsqg} are strong solutions, but 
local existence of strong solutions with general initial data in $\mathcal D(\sqrt{\l})$ is not known. In fact, $\mathcal D(\sqrt{\l})$ might be a borderline case for well posedness. Strong solutions are unique:
\beg{prop}\la{strongunique}
Two solutions of (\ref{sqg}) in $d=2$, obeying (\ref{strongsqg}), and having the same initial data in $\mathcal D(\sqrt{\l})$, must coincide.
\end{prop}
The idea of proof is simple: If $\theta_i$ are the two solutions $i=1,2$,  $\theta$ denotes their difference and $u$ denotes the difference of their velocities, then
\[
\fr{1}{2}\fr{d}{dt}\|\theta\|^2_{L^2} + \|\l^{\fr{1}{2}}\theta\|^2_{L^2} \le 
\int_{\Omega}|u||\theta||\na\Theta|dx
\]
holds, with $\Theta = \fr{1}{2}(\theta_1 + \theta_2)$. The bound
\[
\|\theta(t)\|^2_{L^2} \le \|\theta(0)\|_{L^2}^2\exp{(C\int_0^T\|\nabla \Theta(s)\|^2_{L^2}ds)}
\]
follows because $\mathcal D(\sqrt\l) \subset L^4(\Omega)$ and the Riesz transforms $\na\l^{-1}$ are bounded in $L^4$ (\cite{cabre}).\\
Weak solutions exist globally (\cite{ci}), even without dissipation (\cite{cn}), but are not known to be unique. However, if the initial data are interior Lipschitz continuous, then weak solutions are globally interior Lipschitz continuous. A priori bounds for smooth solutions were given in \cite{sqgb} and a construction was given in \cite{igsqg}.
Let 
\be
d(x) = dist(x,\pa\Omega)
\la{dx}
\ee
denote the distance from $x$ to the boundary of $\Omega$. 

The main result of 
\cite{igsqg} is

\beg{thm}\la{gradint} Let $\theta_0\in H_0^1(\Omega)\cap W^{1,\infty}(\Omega)$ and let $0<T\le\infty$. There exists $\theta(x,t)$, a solution of (\ref{sqg}) on the time interval $[0, T)$, with initial data $\theta(x,0)= \theta_0(x)$ and a constant $\Gamma_1$ depending only on $\Omega$ such that
\be 
\|\theta(\cdot, t)\|_{L^{\infty}(\Omega)} \le \|\theta_0\|_{L^{\infty}(\Omega)},
\la{linftyb}
\ee
and
\be
\sup_{0\le t<T} \sup_{x\in\Omega}d(x)|\na_x\theta(x,t)| \le \Gamma_1\left[\sup_{x\in\Omega}d(x)|\na_x\theta_0(x)| + \left(1+\|\theta_0\|_{L^{\infty}(\Omega)}\right)^{4}\right] := M
\la{gradintb}
\ee
hold.
\end{thm}
This result holds in any dimension $d$. Interior Lipschitz regularity is obtained using nonlinear lower bounds for the square root of the Dirichlet Laplacian (\cite{sqgb}) and commutator estimates. The main obstacle to obtain regularity up to the boundary is the absence of translation invariance, which is most sharply felt near the boundary. The nonlinear lower bounds for the square root of the Dirichlet Laplacian (\cite{ci}, \cite{sqgb}) are similar to those available in the whole space (\cite{cv1}), but have a cut-off due to the boundary. The lack of translation invariance is manifested in the commutator estimates, where the commutator between the square root of the Laplacian and differentiation is of the order $d(x)^{-2}$ pointwise.

In this paper we investigate the behavior of solutions near the boundary. The local solutions obtained in Proposition \ref{locsqg} belong to $C^{\alpha}(\Omega)$ up to the boundary, for any $0<\alpha<1$, but this fact follows from embedding of $H^2(\Omega)\subset C^{\alpha}(\Omega)$ in $d=2$ and the control of the $H^2(\Omega)$ norm is only for short time. An interesting recent work \cite{vas} in the spirit of \cite{caf} shows that a  solution-dependent $C^{\alpha}(\Omega)$ regularity holds as long as the solution is sufficiently smooth. Unfortunately, as we mentioned earlier, smooth solutions can be guaranteed to exist only for a short time. 

The currently available quantitative global in time information for solutions with smooth initial data is comprised of following three components: 

I) Energy bounds, which imply that $\theta\in L^{\infty}(0,T; L^2(\Omega))\cap L^2(0,T; {\mathcal D}(\l^{\fr{1}{2}}))$, 

II) A maximum principle, which implies $\theta\in L^{\infty}(0,T; L^{\infty}(\Omega))$, and,

III) For solutions  constructed by a judicious method mentioned above, the interior Lipschits bound (\ref{gradintb}). \\

No uniqueness is guaranteed. The velocity is given by rotated Riesz transforms. It is known (\cite{cabre}) that if $\theta$ vanishes at the boundary and belongs to $C^{\alpha}$ then its Dirichlet Riesz transforms are in $C^{\alpha}(\Omega)$. If $\theta$ belongs to $C^{\alpha}$ and vanishes at the boundary, then the  stream function $\psi = \l^{-1}\theta$ belongs to $C^{1,\alpha}$ and vanishes at the boundary, and therefore so do its tangential derivatives. Thus, the normal component of the velocity vanishes at the boundary, but no rate is available if $\theta$ belongs to $C^{\alpha}$.

In this work we show that the problem of controlling the H\"{o}lder continuity of the solution up to the boundary depends solely on quantitative bounds on the the vanishing of $\theta$ at the boundary. We prove two results detailing this fact. We consider 
\be
b_1(x,t)  = \fr{\theta(x,t)}{w_1(x)}
\la{b1}
\ee
where $w_1$ is the normalized positive first eigenfunction of the Dirichlet Laplacian, which is known to be smooth and vanish like $d(x)$ at the boundary.
In Theorem \ref{locb1} we show that for solutions constructed from smooth initial data  obeying the a priori information detailed above (I, II, III), and for any $p>d$, there exists a time $T_0$ and a constant $B$, depending only on $\|\theta_0\|_{L^{\infty}}$, $M$ (of \ref{gradintb}) and the initial norm $\|b_1(0)\|_{L^p(\Omega)}$,
such that 
\be
\sup_{0\le t\le T_0}\|b_1(t)\|_{L^p(\Omega)} \le B
\la{b1bound}
\ee
holds. This is a local existence theorem, local because the control of $\|b_1(t)\|_{L^p(\Omega)}$ is maintained for finite time, although the interior Lipschitz bound and the $L^{\infty}$ bound are global. 

Our second main result, Theorem (\ref{glh}), shows that if the bound (\ref{b1bound}) holds for some interval of time, then the solutions constructed in (\cite{igsqg}) are in $C^{\alpha}(\Omega)$ on that interval of time. The H\"{o}lder exponent $\alpha$ is explicit, it is given by $\alpha< 1-\fr{d}{p}$ where $p$ is the exponent in (\ref{b1bound}). Thus, the condition (\ref{b1bound}), which can be maintained for short time, is sufficient for global H\"{o}lder regularity up to the boundary. This condition also implies a quantitative vanishing of the normal component of velocity at the boundary, $u\cdot N  = O(d(x)^{\alpha})$ with rate depending on $M$ and $B$.

The boundedness of $\|b_1\|_{L^p(\Omega)}$ is a weaker condition than $\theta\in W^{1,p_1}_0(\Omega)$, $p_1>p>d$. Thus our condition is necessary for regularity. Local well-posedness in $W^{1,p_1}_0(\Omega)$ is not known. The previously known local existence theory was established in the domain of the Laplacian, which is a strictly smaller space. If the solution is in $W_0^{1,p}(\Omega)$, $p>d$,  then by embedding results, it is in $C^{\alpha}$ up to the boundary for any $0< \alpha <1-\fr{d}{p}$.

In order to prove our results, we obtain key quantitative bounds using $B$.
We show first that if $B$ is finite, then the velocity is bounded (Proposition \ref{uconditional}). By contrast, if only the available apriori information (I, II, III) is used, then the velocity logarithmically diverges with the distance to the boundary (Proposition \ref{ubounduncond}).

Secondly, we obtain bounds for the finite difference quotients of velocity which diverge at the boundary with a sublinear power of the distance, $d(x)^{-\fr{d}{p}}$, (Proposition \ref{deltuh}),  as opposed to $d(x)^{-1}$ in the case of the a priori global information, as was shown in \cite{sqgb}. A quantitative rate of vanishing of the normal component of velocity  is proved in Proposition \ref{unb}.

Thirdly, we obtain bounds for the commutator between finite differences  and $\l$ which diverge subquadratically near the boundary, $d(x)^{-1-\fr{d}{p}}$, (Proposition \ref{commuh})  as opposed to quadratically $d(x)^{-2}$, which is the case in which only the global a priori information (I,II,III) is used. 

These three elements, together with the strong boundary repulsive damping effect of the square root of the Laplacian, form the basis of the proof of persistence of $C^{\alpha}$ regularity, with $\alpha<1-\fr{d}{p}$ .

In the whole space, any $C^{\alpha}$, $\alpha>0$ regularity can be upgraded to 
Lipschitz regularity (and further to $C^{\infty}$ (\cite{cw})). In  bounded domains, while any interior $C^{\alpha}$ regularity can be upgraded to interior Lipschitz regularity (\cite{sqgb}), in general, the problem of global Lipschitz regularity up to the boundary is open.  The passage to Lipschitz bounds up to the boundary is not achievable with our tools, even conditioned on knowledge of linear vanishing of $\theta$ (i.e. even assuming a time-independent bound for $b_1$ in $L^{\infty}$). This is due to the fact that the commutator between derivatives and $\l$ still costs $d(x)^{-1}$ near the boundary.

The paper is organized as follows. After recalling basic facts in Section 
\ref{prel} we prove in in Section \ref{wei} a remarkable generalization of the C\'{o}rdoba-C\'{o}rdoba inequality (\cite{cc}) which was obtained in bounded domains in \cite{ci}. This new pointwise inequality involves weights $w$,
\be
\Phi'(b)\l(wb)-\l(w\Phi(b))\ge (\l(w))\left(b\Phi'(b)-\Phi(b)\right)
\la{remi}
\ee
(see (\ref{wcorcor}, \ref{dphi})) and is valid for any convex function $\Phi$ of one variable which  satisfies $\Phi(0)=0$, any smooth function $b$ and any smooth positive function $w$ which vanishes at $\pa\Omega$. The inequality implies a comparison principle for solutions of drift diffusion equations with Dirichlet square root Laplacian and may have independent interest. We use it with $b = \fr{\theta}{w_1}$ and prove that $B$ of (\ref{b1bound}) persists to be finite if the drift is the sum of a regular function in $L^{\infty}$ whose normal component vanishes at the boundary and a small $L^{\infty}$ function. 
In Section \ref{riesz} we derive bounds for the the Dirichlet Riesz transforms and in Section \ref{gradient} we obtain bounds for finite differences of Dirichlet Riesz transforms. Section \ref{commutator} is devoted to the improved bounds on the commutator between local finite differences and $\l$, and Section \ref{holder} contains the bound for  the H\"{o}lder seminorms near the boundary.

\section{Preliminaries}\la{prel}
We consider $\Omega\subset \Rr^d$ a bounded open set with smooth boundary.
The $L^2(\Omega)$ - normalized eigenfunctions of $-\D$ are denoted $w_j$, and its eigenvalues counted with their multiplicities are denoted $\lambda_j$: 
\be
-\D w_j = \lambda_j w_j.
\la{ef}
\ee
It is well known that $0<\lambda_1\le...\le \lambda_j\to \infty$  and that $-\D$ is a positive selfadjoint operator in $L^2(\Omega)$ with domain ${\mathcal{D}}\left(-\D\right) = H^2(\Omega)\cap H_0^1(\Omega)$.
The ground state $w_1$ is positive and
\be
c_0d(x) \le w_1(x)\le C_0d(x)
\la{phione}
\ee
holds for all $x\in\Omega$, where $c_0, \, C_0$ are positive constants depending on $\Omega$. Functional calculus can be defined using the eigenfunction expansion. In particular
\be
\left(-\D\right)^{\beta}f = \sum_{j=1}^{\infty}\lambda_j^{\beta} f_j w_j
\la{funct}
\ee
with 
\[
f_j =\int_{\Omega}f(y)w_j(y)dy
\]
for $f\in{\mathcal{D}}\left(\left (-\D\right)^{\beta}\right) = \{f\left |\right. \; (\lambda_j^{\beta}f_j)\in \ell^2(\mathbb N)\}$.
We denote by
\be
\l^s = \left(-\D\right)^{\fr{s}{2}}, 
\la{lambdas}
\ee
the fractional powers of the Dirichlet Laplacian, with $0\le s \le 2$
and with $\|f\|_{s,D}$ the norm in ${\mathcal{D}}\left (\l^s\right)$:
\be
\|f\|_{s,D}^2 = \sum_{j=1}^{\infty}\lambda_j^{s}f_j^2.
\la{norms}
\ee
It is well-known that
\[
{\mathcal{D}}\left( \l \right) = H_0^1(\Omega).
\]
Note that in view of the identity
\be
\lambda^{\fr{s}{2}} = c_{s}\int_0^{\infty}(1-e^{-t\lambda})t^{-1-\fr{s}{2}}dt,
\la{lambdalpha}
\ee
with 
\[
1 = c_{s} \int_0^{\infty}(1-e^{-\tau})\tau^{-1-\fr{s}{2}}d\tau,
\]
valid for $0\le s <2$, we have the representation
\be
\left(\left(\l\right)^{s}f\right)(x) = c_{s}\int_0^{\infty}\left[f(x)-e^{t\D}f(x)\right]t^{-1-\fr{s}{2}}dt
\la{rep}
\ee
for $f\in{\mathcal{D}}\left(\left (-\l\right)^{s}\right)$.
We use precise upper and lower bounds for the kernel $H_D(t,x,y)$ of the heat operator,
\be
(e^{t\D}f)(x) = \int_{\Omega}H_D(t,x,y)f(y)dy .
\la{heat}
\ee
These are as follows (\cite{davies1},\cite{qszhang1},\cite{qszhang2}).
There exists a time $T>0$ depending on the domain $\Omega$ and constants
$c$, $C$, $k$, $K$, depending on $T$ and $\Omega $ such that
\be
\ba
c\min\left (\fr{w_1(x)}{|x-y|}, 1\right)\min\left (\fr{w_1(y)}{|x-y|}, 1\right)t^{-\fr{d}{2}}e^{-\fr{|x-y|^2}{kt}}\le \\H_D(t,x,y)\le C
\min\left (\fr{w_1(x)}{|x-y|}, 1\right)\min\left (\fr{w_1(y)}{|x-y|}, 1\right)t^{-\fr{d}{2}}e^{-\fr{|x-y|^2}{Kt}}
\ea
\la{hb}
\ee
holds for all $0\le t\le T$. Moreover
\be
\fr{\left |\na_x H_D(t,x,y)\right|}{H_D(t,x,y)}\le
C\left\{
\ba
\fr{1}{d(x)},\quad\quad \quad\quad {\mbox{if}}\; \sqrt{t}\ge d(x),\\
\fr{1}{\sqrt{t}}\left (1 + \fr{|x-y|}{\sqrt{t}}\right),\;{\mbox{if}}\; \sqrt{t}\le d(x)
\ea
\right.
\la{grbx}
\ee
holds for all $0\le t\le T$.
Note that
\be
H_D(t,x,y) = \sum_{j=1}^{\infty}e^{-t\lambda_j}w_j(x)w_j(y) ,
\la{hphi}
\ee
and therefore long time $t\ge T$ estimates are rather straightforward.
The gradient bounds (\ref{grbx}) result by symmetry in
\be
\fr{\left |\na_y H_D(t,x,y)\right|}{H_D(t,x,y)}\le
C\left\{
\ba
\fr{1}{d(y)},\quad\quad \quad\quad\quad {\mbox{if}}\; \sqrt{t}\ge d(y),\\
\fr{1}{\sqrt{t}}\left (1 + \fr{|x-y|}{\sqrt{t}}\right),\;{\mbox{if}}\; \sqrt{t}\le d(y).
\ea
\right.  
\la{grby}
\ee
We use as well the  bounds (\cite{sqgb})
\be
\na_x\na_x H_D(x,y,t) \le Ct^{-1-\fr{d}{2}}e^{-\fr{|x-y|^2}{\tilde{K}t}}
\la{naxnaxb}
\ee
valid for $t\le cd(x)^2$ and $0<t\le T$,  and
\be
\na_x\na_x H_D(x,y,t) \le Cd(x)^{-2}t^{-\fr{d}{2}}e^{-\fr{|x-y|^2}{\tilde{K}t}}
\la{naxnaxblt}
\ee
for $t\ge cd(x)^2$, which follow from the upper bounds (\ref{hb}), (\ref{grbx}). 
Important additional bounds we use are  
\be
\int_{\Omega}\left |(\na_x +\na_y)H_D(x,y,t)\right|dy \le Ct^{-\fr{1}{2}}e^{-\fr{d(x)^2}{\tilde{K}t}},
\la{cancel1}
\ee
with pointwise version 
\be
\left |(\na_x +\na_y)H_D(x,y,t)\right| \le ct^{-\fr{d+1}{2}}e^{-\fr{d(x)^2}{\tilde{K}t}},
\la{cancel12}
\ee
and
\be
\int_{\Omega}\left |\na_x(\na_x +\na_y)H_D(x,y,t)\right|dy \le Ct^{-1}e^{-\fr{d(x)^2}{\tilde{K}t}},
\la{cancel2}
\ee
with pointwise version 

\be 
\left |\na_x(\na_x +\na_y)H_D(x,y,t)\right| \le Ct^{-\fr{d+2}{2}}e^{-\fr{d(x)^2}{\tilde{K}t}}.
\la{cancel21}
\ee
valid for $t\le cd(x)^2$ and $0<t\le T$. These bounds reflect the fact that translation invariance is remembered in the solution of the heat equation with Dirichlet boundary data for short time, away from the boundary. They were proved in \cite{sqgb}, \cite{cn}.

The following elementary lemma is used in several instances:
\beg{lemma}\la{intpk}
Let $\rho>0$, $p>0$. Then
\be
\int_0^{\rho^2} t^{-1-\fr{m}{2}}\left(\fr{p}{\sqrt{t}}\right)^je^{-\fr{p^2}{Kt}}dt \le  C_{K,m,j}p^{-m}
\la{pbeta}
\ee
if $m\ge 0$, $j\ge 0$, $m+j>0$, and
\be
\int_0^{\rho^2} t^{-1}e^{-\fr{p^2}{Kt}}dt = \int_{\fr{p^2}{K\rho^2}}^{\infty}x^{-1}e^{-x}dx
\la{pzero}
\ee
if $m=0$ and $j=0$, with constants $C_{K,m,j}$ independent of $\rho$ and $p$. Note that when $m+j>0$, $\rho=\infty$ is allowed. Note also that the right-hand side of (\ref{pzero}) is exponentially small if $\rho \le \epsilon p$. 

\end{lemma}

We recall from \cite{ci} that the C\'{o}rdoba-C\'{o}rdoba inequality (\cite{cc}) holds in bounded domains. In fact, more is true: there is a lower bound that provides a strong boundary repulsive term:

\beg{prop}{\la{cordoba}} Let $\Omega$ be a bounded domain with smooth boundary. Let $0\le s<2$. There exists a constant $c>0$ depending only on the domain $\Omega$ and on $s$, such that, for any
$\Phi$, a $C^2$ convex function satisfying $\Phi(0)= 0$, and any  $f\in C_0^{\infty}(\Omega)$, the inequality
\be
\Phi'(f)\l^s f - \l^s(\Phi(f))\ge \fr{c}{d(x)^s}\left(f\Phi'(f)-\Phi(f)\right)
\la{cor}
 \ee
holds pointwise in $\Omega$.
\end{prop}
We specialize from now on to $s=1$. We use in particular the result above in the form (\cite{ci})  
\be
D(f)(x) = \left(f\l f - \fr{1}{2}\l\left({f^2}\right)\right)(x) \ge \gamma_1\fr{f^2(x)}{d(x)}
\la{dfdxb}
\ee
with $\gamma_1>0$ depending only on $\Omega$.

\section{Weighted estimates}\la{wei} 
Let $w(x)$ be a function which is positive in $\Omega$ and belongs to $\mathcal D (\l)$, for instance $w(x) = w_1(x)$.
\beg{lemma}\la{wcordoba} Let $\Phi$ be a convex function of one variable with $\Phi(0) = 0$. Then
\be
\Phi'(b(x))\l(wb)(x) - \l(w\Phi(b))(x)  = (\l w)(x)\left(b(x)\Phi'(b(x)) -\Phi(b(x)\right) + D_\Phi(x)
\la{wcorcor}
\ee
with
\be
D_\Phi(x) = c\int_0^{\infty}t^{-\fr{3}{2}}\int_{\Omega}w(y)H_D(x,y,t)\left[\Phi(b(y))-\Phi(b(x)) - \Phi'(b(x))(b(y)-b(x))\right]dydt
\la{dphi}
\ee
\end{lemma}

\noindent{\bf{Proof.}} Let $\phi(x) = \Phi(b(x))$ and $\phi'(x) = \Phi'(b(x))$.We have
\[
\ba
(\phi'\l(wb)- \l(w\phi))(x) =
c\int_0^{\infty}t^{-\fr{3}{2}}\left[\phi'(x)w(x)b(x) - \int_{\Omega}\phi'(x)w(y)H_D(x,y,t)b(y)dy\right ]dt\\
- c\int_0^{\infty}t^{-\fr{3}{2}}\left[w(x)\phi(x) - \int_{\Omega}w(y)H_D(x,y,t)\phi(y)dy\right]dt\\
= (b(x)\phi'(x)-\phi(x))c\int_0^{\infty}t^{-\fr{3}{2}}\left[w(x)- \int_{\Omega}w(y)H_D(x,y,t)dy\right]dt + D_\Phi(x)\\=
(b(x)\phi'(x)-\phi(x))\l(w)(x) + D_{\Phi}(x)
\ea
\]
\beg{rem}
We note that $D_{\Phi}\ge 0$ for convex functions $\Phi$ because the integrand is non-negative.
\end{rem}
Let us consider now an evolution equation
\be
\pa_t\theta + v\cdot\na \theta + \l\theta = 0
\la{vthetaeq}
\ee
with $v = v(x,t)$ a divergence-free vector field tangent to the boundary of $\Omega$. Let us consider a smooth enough weight $w(x,t)>0$ which vanishes at the boundary of $\Omega$  and compute the evolution of $\Phi(b(x,t))$
where $\Phi$ is a nonnegative convex function of one variable, with $\Phi(0)=0$, and where
\be
b(x,t) = \fr{\theta(x,t)}{w(x,t)}.
\la{b}
\ee
In view of (\ref{wcorcor}), we obtain the remarkable equation
\be
(\pa_t + v\cdot\na + \l)(w\Phi(b)) + \left((\pa_t +v\cdot\na +\l)w\right)(b\Phi'(b)-\Phi(b)) + D_{\Phi} = 0,
\la{remeq}
\ee
where $D_{\Phi}$ is defined above in (\ref{dphi}). Denoting 
\be
L_v = \pa_t + v\cdot\na + \l
\la{lv}
\ee
we have thus
\be
L_v(w\Phi(b)) + (L_v(w))(b\Phi'(b)- \Phi(b)) + D_{\Phi} = 0
\la{lvineq}
\ee
for $w>0$ and $\Phi$ convex with $\Phi(0)=0$.
There are several important consequences of this identity. In view of the fact that
\be
\int_{\Omega}\l(w\Phi(b))dx = \int_{\Omega}w\Phi(b)\l(1)dx
\la{int1}
\ee
and the lower bound
\be
(\l 1)(x) \ge c_0\fr{1}{w_1(x)}
\la{lowl1}
\ee
we have that
\be
\int_{\Omega}\l(w\Phi(b))dx \ge c_0\int_{\Omega}\left(\fr{w(x)}{w_1(x)}\right) \Phi(b(x))dx. 
\la{lowlwphi}
\ee
Therefore 
\be
\fr{d}{dt}\int_{\Omega}w\Phi(b)dx + c_0\int_{\Omega}\left(\fr{w(x)}{w_1(x)}\right) \Phi(b(x))dx  +\int_{\Omega}(L_v(w))(b\Phi'(b)-\Phi(b))dx + \int_{\Omega} D_{\Phi}(x,t)dx \le 0.
\la{phlvineq}
\ee
Let us take now $\Phi$ to be (a smooth convex approximation of) the function
\be
\Phi_B(b) = (b-B)_{+}
\la{PhiB}
\ee
where $B$ is a large fixed number. Notice that in this case
\be
b\Phi_B'(b)-\Phi_B(b) = BH(b-B),
\la{vanishae}
\ee
where $H(x)$ is the Heaviside function. Because
$b\Phi'(b)-\Phi(b)\ge 0$, if $L_v(w)\ge 0$, then 
\be
\fr{d}{dt}\int_\Omega w(x,t)\Phi_B(b(x,t))dx \le 0.
\la{phiBintdecay}
\ee
It follows that, If $\Phi(b(x,0)) = 0$, then $\Phi(b(x,t)) = 0$ for $t\ge 0$.
Applying this reasoning to the functions $b$ defined above in (\ref{b}) as well as to $b_{-} = \fr{-\theta}{w_1}$ we obtain
\be
|\theta(x,t)|\le Bw(x,t).
\la{thetawbound}
\ee
\beg{rem} 
This shows that if $L_v(w) \ge 0$ and $|\theta_0(x)|\le Bw(x,0)$, then (\ref{thetawbound}) holds.
\end{rem}
\beg{thm}\la{w1bound} Let $\theta$ solve (\ref{vthetaeq}) where $v$ is a continuous, divergence-free field, tangent to the boundary. 
Assume that there exists a constant $\gamma(t)$ such that
\be
v\cdot\na w_1 + \gamma(t)w_1\ge 0
\la{toomuch}
\ee
holds for $x\in\Omega$ and $t\ge 0$. Assume that the initial data $\theta_0$ obeys
\be
|\theta_0(x)| \le Bw_1(x)
\la{thethidB}
\ee
for all $x\in\Omega$.
Then
\be
|\theta(x,t)| \le Bw_1(x)e^{-t\sqrt{\lambda_1} + \int_0^t\gamma(s)ds}
\la{thetaB}
\ee
holds for all $x\in\Omega$ and all $t\ge 0$.
\end{thm}
\noindent{\bf{Proof.}} Consider 
\be 
w(x,t) = e^{-t\sqrt{\lambda_1} + \int_0^t\gamma(s)ds}w_1.
\la{wt}
\ee
Note that the assumption (\ref{toomuch}) implies that
\be
L_v(w(x,t)) \ge 0
\la{lvwgez}
\ee
Then we use (\ref{thetawbound}) and conclude the proof.

If $v$ is bounded and if its normal component vanishes of first order
at the boundary then the condition (\ref{toomuch}) is satisfied.

\beg{prop} \la{sufcon} Condition (\ref{toomuch}) is satisfied if
$v$ is bounded,
\be
\|v(t)\|_{L^{\infty}}\le V(t),
\la{vinfty}
\ee
and has a normal component which vanishes to first order near the boundary of $\Omega$,
\be 
|v(x,t)\cdot N(x)|  \le V(t)d(x),
\la{vanishvN}
\ee
where  $N(x)$ is a continuous unit vector defined near the boundary $\pa\Omega$ and extending the normal at $\pa\Omega$. 
\end{prop}
 Indeed, for any smooth vector field $T(x)$ defined near the boundary and tangent to the boundary,  we have by the smoothness of $w_1$ and its equivalence to the distance to the boundary,
\be
T(x)\cdot\na w_1(x) \le Cw_1(x)
\la{tw1}
\ee
near the boundary $\Omega$. Then we decompose $v= (v\cdot T)T + (v\cdot N)N = v_T + v_N$ with $T$ smooth near the boundary and tangent to the boundary, and use the fact that
\be
|v_N(x)|\le Cw_1(x),
\la{vNb}
\ee
near the boundary, which follows by the assumption (\ref{vanishvN}). The fact that $|v\cdot\na w_1|\le C w_1$ away from the boundary follows from the boundedness of $v$. This concludes the proof of Proposition \ref{sufcon}. 
\beg{rem}
Theorem \ref{w1bound} can be proved also using
\be
\Phi(b) = b^{2m}.
\la{powerphi}
\ee
We note that in this case
\be
b\Phi'(b)-\Phi(b) = (2m-1)\Phi(b).
\la{phipm}
\ee
We take $w=w_1$ and use the fact that
\be
L_v(w_1) = v\cdot\na w_1 + \sqrt{\lambda_1}w_1,
\la{lvw1}
\ee
and returning to (\ref{phlvineq}) we obtain 
\be
\fr{d}{dt}\int_{\Omega} w_1(x)\Phi(b(x,t))dx \le (2m-1)(\gamma(t) -\sqrt{\lambda_1})\int_{\Omega}w_1(x)\Phi(b(x,t))dx
\la{powgammaineq}
\ee
where
\be
\gamma(t) = \sup_{x\in\Omega}\left (-\fr{v(x,t)\cdot\na w_1(x)}{w_1(x)}\right)
\la{gammat}
\ee
Integrating in time, taking $2m$ roots and then the limit $m\to\infty$, we arrive at
\be
\|b(t)\|_{L^{\infty}}\le \|b_0\|_{L^{\infty}}e^{- t\sqrt{\lambda_1} + \int_0^t\gamma(s)ds}.
\la{btgammabound}
\ee

Note that if (\ref{toomuch}) holds then (\ref{btgammabound}) is precisely (\ref{thetaB}).
\end{rem}
We consider now the case of fixed $m$. 
\beg{prop}\la{fixedm} Let $m\ge 1$ be an integer, let $v$ be a bounded divergence-free function which can be decomposed
\be
v= v_r + v_s
\la{vrs}
\ee
with $v_r(x,t)$ obeying $\gamma_r\in L^1[0,T]$, where $\gamma_r(t)$ is defined as in (\ref{gammat}) by
\be
 \sup_{x\in\Omega}\left (-\fr{v_r(x,t)\cdot\na w_1(x)}{w_1(x)}\right) = \gamma_r(t)
\la{gammart}
\ee
and with 
\be
\|v_s(t)\|_{L^{\infty}} \le \fr{c_0}{(2m-1)\|\na w_1\|_{L^{\infty}}}
\la{vsbound}
\ee
where $c_0$ is the constant from (\ref{phlvineq}). Then
\be
\int_{\Omega} w_1(x) \left(\fr{\theta(x,t)}{w_1(x)}\right)^{2m}dx \le 
  e^{(2m-1)\left(-t\sqrt{\lambda_1} +\int_0^t\gamma_r(s)ds\right)}\int_{\Omega}w_1(x) \left(\fr{\theta_0(x)}{w_1(x)}\right)^{2m}dx 
\la{bmbound}
\ee
holds for $t\in [0,T]$.
\end{prop}
\beg{rem} Note that the right hand side of (\ref{vsbound}) depends only on $\Omega$ and $m$.
\end{rem}

\noindent{\bf Proof.} The proof follows along the same lines as above. We take 
$\Phi$ as in (\ref{powerphi}), 
$w=w_1$, and using the decomposition we have that
\be
L_v(w_1) = v_s\cdot\na w_1 + \left(\sqrt{\lambda_1} + \left(\fr{v_r\cdot\na w_1}{w_1}\right)\right )w_1.
\la{lvrsw1}
\ee
Consequently, from (\ref{gammart}) and (\ref{vsbound}) we have
\be
(2m-1)L_v(w_1) \ge -c_0 + (2m-1)(\sqrt{\lambda_1} -\gamma_r(t))w_1(x).
\la{lvrsw1bound}
\ee
We use this inequality and (\ref{phipm}) in (\ref{lvineq}), integrate in time, and deduce (\ref{bmbound}).

We record here a lemma relating weighted and unweighted norms of $b$:
\beg{lemma}\la{weightnotweight} Let $m>p\ge 1$. Then, there exists a constant $C_{m,p}$ depending only on $\Omega$, $m$ and $p$ such that 
\be
\|b\|_{L^p(\Omega)} \le C_{m,p}\left(\int_{\Omega}w_1(x)b^{2m}(x)dx\right)^{\fr{1}{2m}}
\la{bwb}
\ee
holds for any $b$. Conversely, let $p\ge 2m-1\ge 1$ and let $b_1 =\fr{\theta}{w_1}$. Then 
\be
\left(\int_{\Omega}w_1(x)b_1^{2m}(x)dx\right)^{\fr{1}{2m}}
\le \|\theta\|_{L^{\infty}(\Omega)}^{\fr{1}{2m}}\|b_1\|_{L^{p}(\Omega)}^{\fr{2m-1}{m}}|\Omega|^{\fr{p+1-2m}{2mp}}
\la{bpb}
\ee
\end{lemma}
\noindent{\bf{Proof.}} The first inequality uses just the H\"{o}lder inequality for the functions  $w_1(x)^{\fr{p}{2m}}|b(x)|^p$ and $w_1(x)^{-\fr{p}{2m}}$, with exponents
$\fr{2m}{p}$, $\fr{2m}{2m-p}$, and 
\be
A_{m,p} =\int_{\Omega} w_1(x)^{-\fr{p}{2m-p}}dx <\infty
\la{w1negint}
\ee
which holds because  $\fr{p}{2m-p}<1$. Then $C_{m,p} = A_{m,p}^{\fr{2m-p}{2mp}}$. The second inequality is straighforward.

\section{Bounds for Riesz transforms}\la{riesz}
We consider $u$ given in (\ref{u}),  
\[
u = \na^{\perp}\l^{-1}\theta.
\]
We are interested in estimates of $u$ in terms of $\theta$.

\beg{prop}\la{ubounduncond} Let $u$ be given by (\ref{u}) and let $\theta$ be bounded and interior Lipschitz, i.e., obeying 
\be
d(y)|\na\theta(y)| \le M.
\la{M}
\ee
Then, there exist  constants $C$ depending only on the domain $\Omega$ such that
\be
|u(x)| \le CM + C\|\theta\|_{L^{\infty}}\left(1 + \log\left(\fr{C}{d(x)}\right)\right).
\la{uuncondb}
\ee
As a consequence, there exist constants $\gamma>0$ and $C$, depending only on the domain $\Omega$, $M$ and $\|\theta\|_{L^{\infty}}$ such that
\be
\int_{\Omega}e^{\gamma |u(x)|}dx \le C.
\la{expubound}
\ee
\end{prop}
\beg{rem}
The bound (\ref{uuncondb}) does not use any information about vanishing of $\theta$ at the boundary, but it uses (\ref{M}) which follows in our case from a priori bounds (\ref{gradintb}). The bound is in particular true for $\theta =1$, where we know that the Riesz transform is in general only BMO and is not bounded all the way to the boundary.
\end{rem}

\noindent{\bf Proof.} In view of the representation
\be
\l^{-1} = c\int_0^{\infty}t^{-\fr{1}{2}}e^{t\D}dt
\la{lambdaminusone}
\ee
we have that
\be 
u(x) = c\int_0^{\infty}t^{-\fr{1}{2}}dt\int_{\Omega}\na_x^{\perp}H_D(x,y,t)\theta(y)dy.
\la{uintegr}
\ee
We split
\be
u = u^{in} + u^{out}
\la{usplit}
\ee
with
\be
u^{in}(x) = c\int_0^{\rho^2}t^{-\fr{1}{2}}dt\int_{\Omega}\na_x^{\perp}H_D(x,y,t)\theta(y)dy
\la{uint}
\ee
with 
$\rho=\rho(x)$ a length scale  smaller than the distance from $x$ to the boundary of $\Omega$:
\be
\rho \le \epsilon d(x).
\la{rholess}
\ee
We split further
\be
\ba
u^{in}(x) = 
\int_0^{\rho^2}t^{-\fr{1}{2}}\int_{\Omega}(\na_x^{\perp} + \na_y^{\perp})H_D(x,y,t)\theta(y)dydt  + v(x)\\
= u^{in}_1(x) + v(x)
\ea
\la{uinv}
\ee 
and then
\be
\ba
u^{in}_1(x) = \int_0^{\rho^2}t^{-\fr{1}{2}}\int_{\Omega}(\na_x^{\perp} + \na_y^{\perp})H_D(x,y,t)(\phi (y) + (1-\phi(y))) \theta(y)dydt\\
= v_1 +v_2
\ea
\la{vs}
\ee 
where $\phi$ is a standard cutoff compactly  supported in a ball of radius $\ell = \fr{\eta d(x)}{4}$ around $x$ and identically one in the ball of radius $\fr{\ell}{4}$. Above $\epsilon>0$ and $\eta>0$ are small numbers at our disposal.

We use (\ref{cancel12}) to bound
\[
t^{-\fr{1}{2}}\left|(\na_x^{\perp}+ \na_y^{\perp})H_D(x,y,t)\right |\le C(d(x))^{-(d+2)}
\]
and thus
\be
|v_1(x)| = \left |\int_0^{\rho^2}t^{-\fr{1}{2}}\int_{\Omega}(\na_x^{\perp} + \na_y^{\perp})H_D(x,y,t)\phi(y)\theta(y)dy\right| \le C\epsilon^2\eta^d\|\theta\|_{L^{\infty}}.
\la{v1bound}
\ee
Here we estimate the volume of the support of $\phi$ by $C\eta^dd(x)^d$.
For $v_2$ we use the bounds (\ref{grbx}) and (\ref{grby}) and the fact that $|x-y|\ge \eta d(x)/8$ on the support of $1-\phi$ to obtain
\be
|v_2(x)| \le C\frac{\epsilon^2}{\eta^2}\|\theta\|_{L^{\infty}}.
\la{v2bound}
\ee
Here we used the fact that $t^{-\fr{d+2}{2}}e^{-\fr{|x-y|^2}{Kt}}\le C|x-y|^{-(d+2)}$, and
\[
\rho^2\int_{|x-y|\ge \eta d(x)/8}|x-y|^{-(d+2)} \le C\frac{\epsilon^2}{\eta^2}.
\]
Now we write
\be
v(x) = -c\int_0^{\rho^2}t^{-\fr{1}{2}}dt\int_\Omega \na^{\perp}_y H_D(x,y,t)(\phi(y) + 1-\phi(y))\theta(y)dy = v_3+v_4
\la{v34}
\ee
We observe that $v_4$ is estimated exactly like $v_2$, 
\be
|v_4(x)| \le C\frac{\epsilon^2}{\eta^2}\|\theta\|_{L^{\infty}}
\la{v4bound}
\ee
In $v_3$ we integrate by parts
\be
v_3(x) = c\int_0^{\rho^2}t^{-\fr{1}{2}}dt\int_{\Omega}H_D(x,y,t)\na^{\perp}(\phi(y)\theta(y))dy. 
\la{v3}
\ee
We use here $|\na\phi|\le C(\eta d(x))^{-1}$, and (\ref{M}). Because on the support of $\phi$ we have $d(x) \le 2d(y)$,
we deduce that
\be
|v_3(x)|\le C\fr{M+\eta^{-1}\|\theta\|_{L^{\infty}}}{d(x)}\int_0^{\rho^2(x)}t^{-\fr{1}{2}}dt \le C\epsilon(M + \eta^{-1}\|\theta\|_{L^{\infty}}),
\la{v3bound}
\ee
and consequently, from (\ref{v1bound}), (\ref{v2bound}), (\ref{v4bound}) and 
(\ref{v3bound}) that
\be
|u^{in}(x)|\le C\left(\epsilon M + (\fr{\epsilon^2}{\eta^2} + \epsilon^2\eta^d + \fr{\epsilon}{\eta})\|\theta\|_{L^{\infty}}\right)
\la{uintbound}
\ee
We write  $u^{out}$  as
\be
u^{out}(x) = u_T(x) + u^{T}(x)
\la{uouT}
\ee
where
\be
u_T(x) = \int_{T}^{\infty} t^{-\fr{1}{2}}\int_{\Omega}\na_{x}^{\perp}H_D(x,y,t)\theta(y)dy
\la{ufromT}
\ee
and
\be
u^{T}(x) = \int_{\rho^2(x)}^T t^{-\fr{1}{2}}\int_{\Omega}\na_{x}^{\perp}H_D(x,y,t)\theta(y)dy.
\la{uuptoT}
\ee
Because (\ref{hphi}) $u_T$ is smooth, and in particular 
\be
\|u_T\|_{L^{\infty}}\le C\|\theta\|_{L^{\infty}}
\la{ufromTb}
\ee
We take $\epsilon=\eta =1$. From (\ref{uintbound}) and (\ref{ufromTb}) we have
\be
|u^{in}(x)| + |u_T(x)| \le C(\|\theta\|_{L^{\infty}} + M).
\la{uintTb}
\ee
On the other hand, from (\ref{phione}) and (\ref{hb}) we bound
\be
\fr{1}{d(x)} H_D(x,y,t) \le C\fr{1}{|x-y|}t^{-\fr{d}{2}}e^{\fr{|x-y|^2}{Kt}}
\la{hddxb}
\ee
and using (\ref{grbx}) we obtain
\be
\ba
|u^T(x)| \le c\int_{\rho^2}^T t^{-\fr{1}{2}}dt\int_{\Omega}t^{-\fr{d}{2}}e^{-\fr{|x-y|^2}{Kt}}\fr{1}{|x-y|}|\theta(y)|dy\\
\le C\|\theta\|_{L^{\infty}}\left(\int_{\Rr^2}|y|^{-1}e^{-|y|^2}dy\right)\int_{\rho^2}^Tt^{-1}dt\\
= C\log\left(\fr{T}{\rho^2(x)}\right)\|\theta\|_{L^{\infty}}.
\ea
\la{uuptoTb}
\ee

 For this result we  used  $\epsilon=\eta =1$, but $\epsilon$  may be used to show that the dependence on $M$ is logarithmic for large M. This concludes the proof of Proposition \ref{ubounduncond}.

The next proposition uses information about vanishing of $\theta$ at the boundary.
\beg{prop}\la{uconditional}
Let  $\theta$ be bounded and interior Lipschitz, i.e., obeying (\ref{M}). Let
\be
b_1(x) = \fr{\theta(x)}{w_1(x)}.
\la{bone}
\ee
Let $u$ be given by (\ref{u}). For any $p> d$, there exist constants $C$ depending only on the domain $\Omega$ and $p$ such that
\be
\|u\|_{L^{\infty}}\le CM + C\|\theta\|_{L^{\infty}}\left(1 + \log \|b_1\|_{L^p(\Omega)}\right).
\la{ucondb}
\ee
\end{prop}

\noindent{\bf Proof}. We proceed like in the proof of Proposition \ref{ubounduncond}, and in particular we use the bound (\ref{uintTb}). We bound $u^T$ differently. We take a small number $\delta$ and we split
\be
u^T = u_1 + u_2
\la{uTonetwo}
\ee
where 
\be 
u_1(x) = \int_{\rho^2}^T t^{-\fr{1}{2}}\int_{\Omega\cap\;|x-y|\le \delta}\na_x^{\perp}H_D(x,y,t)\theta(y)dy
\la{uoneout}
\ee
and
\be
u_2(x) = \int_{\rho^2}^T t^{-\fr{1}{2}}\int_{\Omega\cap\;|x-y|\ge \delta}\na_x^{\perp}H_D(x, y,t)\theta(y)dy
\la{utwoout}
\ee
For $u_2$  we have, using (\ref{grbx}), the bound (\ref{hddxb})
and Lemma \ref{intpk} with $j=0$ and $m =d-1$,
\be
|u_2(x)| \le C\|\theta\|_{L^{\infty}}\int_{\Omega \cap \;|x-y|\ge \delta}|x-y|^{-{d}}e^{-\fr{|x-y|^2}{2KT}}dy \le C\|\theta\|_{L^{\infty}}e^{-\fr{\delta^2}{2KT}}\log\left(\fr{C}{\delta}\right).
\la{u2outbound}
\ee 
In order to estimate $u_1$ we write
\be
\theta (y) = b_1(y)w_1(y)
\la{bonetheta}
\ee
and use
\be
w_1(y) \le Cd(y) \le C( d(x) + |x-y|)
\la{triangle}
\ee
and (\ref{hddxb}) to estimate 
\be
\ba
|u_1(x)| \le C\int_{\rho^2}^T t^{-\fr{1}{2}}\int_{\Omega\cap\;|x-y|\le \delta}|b_1(y)|\left(1 +\fr{d(x)}{|x-y|}\right)t^{-\fr{d}{2}}e^{-\fr{|x-y|^2}{Kt}}dy\\
\le C\int_{\Omega\cap\;|x-y|\le \delta}|b_1(y)||x-y|^{-(d-1)}dy + C\fr{1}{\epsilon}\int_{\rho^2}^T t^{-\fr{d}{2}}e^{-\fr{|x-y|^2}{Kt}}dt\int_{\Omega\cap\;|x-y|\le \delta}|b_1(y)||x-y|^{-1}dy\\
\le C(1+\fr{1}{\epsilon})\int_{\Omega\cap\;|x-y|\le \delta}|b_1(y)||x-y|^{-(d-1)}dy 
\ea
\la{u1outbound}
\ee
where we used  Lemma \ref{intpk} with $j=0$ and $m = d-1$ in the first term of the second inequality and $d(x)\le \epsilon^{-1}\sqrt{t}$ and Lemma \ref{intpk} with $j=0$ and $m = d-2$ in the second term. If $d=2$ we treat the second term of the second inequality directly, ignoring the exponential and integrating 
\[
d(x)\int_{\rho^2}^T t^{-\fr{3}{2}}dt \le \fr{2}{\epsilon}.
\]

From the bounds above we obtain
\be
|u^T(x)| \le C\|\theta\|_{L^{\infty}}e^{-\fr{\delta^2}{2KT}}\log\left(\fr{C}{\delta}\right) + C\left(1+\fr{1}{\epsilon}\right)\int_{\Omega\cap\;|x-y|\le \delta}|b_1(y)||x-y|^{-(d-1)}dy.
\la{uuptoTbounde}
\ee

We take now $\eta=\epsilon=1$.
Putting together the estimates  (\ref{uintTb}) and (\ref{uuptoTbounde}) we obtain 

\be
|u(x)| \le C(M +\|\theta\|_{L^{\infty}}) + C\|\theta\|_{L^{\infty}}\log\left(\fr{C}{\delta}\right) + C\int_{\Omega\cap\;|x-y|\le \delta}|b_1(y)||x-y|^{-(d-1)}dy.
\la{udeltabound}
\ee
The estimate (\ref{ucondb}) follows by appropriately choosing $\delta$ small enough. This ends the proof of Proposition \ref{uconditional}. Clearly, the condition $b_1\in L^p(\Omega)$, $p>d$ can be relaxed to
\be
\lim_{\delta\to 0}\int_{\Omega\cap |x-y|\le \delta}|b_1(y)||x-y|^{-(d-1)}dy = 0.
\la{cond}
\ee
We show now a decomposition of the type (\ref{vrs}).

\beg{prop}\la{decomprop} If $p>d$ and if $b_1 = \fr{\theta}{w_1}\in L^p(\Omega)$ then, for any $c_r>0$, there exists $\tau>0$ depending only on $\Omega$, the norm $\|b_1\|_{L^p(\Omega)}$, the constant 
$M$ of (\ref{M}) and on $\|\theta\|_{L^{\infty}}$ such that
\be
u_s(x) := \int_0^{\tau} t^{-\fr{1}{2}}dt\int_\Omega \na_x^{\perp}H_D(x,y,t)\theta(y)dy
\la{us}
\ee
obeys
\be
\|u_s\|_{L^{\infty}} \le c_r.
\la{ussmall}
\ee
\end{prop}
\noindent{\bf{Proof}}. 
We note that
\be 
u_s(x) = u^{in}(x) + u^{\tau}(x)
\la{usintau}
\ee
where $u^{in}$ is defined in (\ref{uint}) and hence obeys (\ref{uintbound}) and where $u^{\tau}$ is $u^T$ of (\ref{uuptoT}) for $T=\tau$. We choose $\eta>0$ of order one, then $\epsilon$ to be small enough such that $\epsilon \eta^{-1}$ also is small enough, so that from (\ref{uintbound}) we obtain
\be
|u^{in}(x)|\le \fr{c_r}{4}.
\la{uincr}
\ee
With these choices of $\epsilon$ and $\eta$ we use the fact that 
\be
\int_{\Omega\cap |x-y|\le \delta}|b_1(y)||x-y|^{-(d-1)}dy \le \|b_1\|_{L^p(\Omega)}\delta^{1-\fr{d}{p}}
\la{deltap}
\ee
to choose $\delta$ small enough so that
\be
C\left(1+\fr{1}{\epsilon}\right)\int_{\Omega \cap\;|x-y|\le \delta}|b_1(y)||x-y|^{-1}dy\le \fr{c_r}{4}.
\la{bdeltasmall}
\ee
Now, once these choices have been made , we choose $\tau$ small enough so that
\be
 C\|\theta\|_{L^{\infty}}e^{-\fr{\delta^2}{2K\tau}}\log\left(\fr{C}{\delta}\right) \le \fr{c_r}{2}.
 \la{etachoice}
\ee
The result then follows from the bounds (\ref{uincr}), (\ref{bdeltasmall}), (\ref{etachoice}) and (\ref{uuptoTbounde}). This concludes the proof of Proposition \ref{decomprop}.

We state now the result of local control of $\|b_1\|_{L^p}$.

\beg{thm}\la{locb1} Let $\theta_0\in H_0^{1}(\Omega)\cap W^{1,\infty}(\Omega)$. Let $m>d$. There exists a time $T_0$ depending only on $\|\theta_0\|_{L^{\infty}}$, $\sup_{x\in\Omega}d(x)|\na\theta_0(x)|$ and $\|b_1(0)\|_{L^{2m}(w_1dx)}$ and a solution $\theta(x,t)$ of (\ref{sqg}) obeying (\ref{gradintb}) and 
\be
\|b_1(t)\|_{L^{2m}(w_1dx)}= \left(\int_{\Omega}w_1(x)b_1(x,t)^{2m}dx\right)^{\fr{1}{2m}} \le C 
\la{wb1mbound}
\ee
for $t\le T_0$. Consequently, for $d<p<m$ there exists $B_p$ such that
\be
\sup_{0\le t\le T_0}\|b_1(t)\|_{L^p(\Omega)}\le B_p
\la{b1Bshort}
\ee
holds.
\end{thm} 
\noindent{\bf{Proof}}. The Proposition \ref{fixedm} and Proposition \ref{decomprop} are used in conjunction with Theorem \ref{gradint} and Lemma \ref{weightnotweight}.

\section{Bounds for finite differences of Riesz transforms}\la{gradient}
We consider now finite differences
\be
\delta_h u(x) = c\int_0^{\infty}t^{-\fr{1}{2}}dt\int_{\Omega}\delta_h^x\na_x^{\perp}H_D(x,y,t)\theta(y)dy
\la{duh}
\ee
with $|h| \le \fr{d(x)}{16}$.

\beg{defi}  
Let us consider a small length $\ell_0$, and take $0\le \ell\le \ell_0$. We consider a ball $B$ centered at a point $x_0$ with $d(x_0)\ge 2\ell$ and of radius
$\ell$. We take a smooth nonnegative function $\phi = \Psi\left(\fr{|x-x_0|}{\ell}\right)$ with $\Psi$ a smooth, nonincreasing function of $z\in \Rr_+$, $\Psi(z) =1$ for $z\le \fr{5}{16}$ and $\Psi(z) = 0$ for $z \ge \fr{7}{16}$. We also take  a function $\chi = \Psi\left (\fr{|x-x_0|}{2\ell}\right )$, noting that $0\le \phi\le \chi\le 1$, $\chi \phi = \phi$ and that the support of $1-\chi$ is included in $|x-x_0|\ge \fr{5\ell}{8}$, and  that of $\na \chi$  in $\fr{5\ell}{8}\le |x-x_0| \le \fr {7\ell}{8}$, so they both are disjoint from the support of $\phi$ which is included in $|x-x_0|\le \fr{7\ell}{16}$.
We refer to $\phi$ as a ``standard cutoff with scale $\ell$'' and center $x_0$, and to $\chi$ as its ``companion''.
\end{defi}

\beg{prop}\la{deltuh}
Let $\phi$ be a standard cutoff with scale $\ell$ and companion $\chi$, and let $u$ be given by (\ref{u}). Then for any $\epsilon>0$, there exists $\delta(\epsilon)$ with $\lim_{\epsilon\to 0}\delta(\epsilon) = 0$, and a constant $C_\epsilon$ depending only on $\Omega$, $\epsilon$  and $p$ such that
 \be
\ba
|\phi(x)\delta_h u(x)| \\
\le \sqrt{\epsilon d(x) D(\chi\delta_h\theta)} + C_{\epsilon}|h|d(x)^{-\fr{d}{p}}\|b_1\|_{L^{p}(\Omega)} + \delta(\epsilon) \phi(x)|\delta_h\theta(x)|
\ea
\la{deltauhb}
\ee
holds pointwise.
\end{prop}
Above $D(\chi\delta_h\theta)$ is given in (\ref{dfdxb}).

\noindent{\bf{Proof.}} We start with bouunds for the gradient. We use the representation
\be
\na u(x) = \na u^{in}(x) + \na u^{out}(x)
\la{nauinout}
\ee
where
\be
\na u^{in}(x) = c\int_0^{\rho^2}t^{-\fr{1}{2}}dt\int_{\Omega}\na_x\na_x^{\perp}H_D(x,y,t)\theta(y)dy
\la{nauin}
\ee
and $\rho= \rho(x) \le \epsilon d(x)$. For $\na u^{out}$ we split in three parts. 
The inner portion of the integral 
\[
\int_{\rho^2}^{\epsilon d(x)^2}t^{-\fr{1}{2}}dt\int_{|x-y|\le d(x)}\na_x\na_x^{\perp}H_D(x,y,t)\theta(y)dy
\]
is bounded using (\ref{naxnaxb}) and ignoring the exponential, yielding
\be
\ba
\int_{\rho^2}^{\epsilon d(x)^2}t^{-\fr{d+3}{2}}dt\int_{|x-y|\le d(x)}|b_1(y)|(d(x) + |x-y|)dy\\
\le C\rho^{-(d+1)} d(x)\int_{|x-y|\le d(x)}|b_1(y)|dy \le Cd(x)^{-\fr{d}{p}}\left(\fr{d(x)}{\rho}\right)^{d+1}\|b_1\|_{L^p}.
\ea
\la{nauoutininb}
\ee
We use $|\na_x\na_x H_D (x,y,t)| \le Ct^{-\fr{d}{2}}d(x)^{-2}$ for $t\ge cd^2(x)$ and $|x-y|\le d(x)$
to bound the  integral
\be
\int_{\epsilon d(x)^2}^{\infty}t^{-\fr{d+1}{2}}dt\int_{|x-y|\le d(x)}d(x)^{-2}|b_1(y)|(d(x) + |x-y|)dy\le C_{\epsilon}d(x)^{-\fr{d}{p}}\|b_1\|_{L^p}.
\la{nauoutinb}
\ee
From  $|\na_x\na_x H_D(x,y,t)| \le Cd(x)^{-2} t^{-\fr{d}{2}} e^{-\fr{|x-y|^2}{Kt}}$ for $t\ge Cd(x)^2$, $|x-y|\ge d(x)$ and Lemma \ref{intpk} with $m=d-1$ we obtain
\be
\ba
\int_{0}^{\infty}t^{-\fr{d+1}{2}}\int_{|x-y|\ge d(x)}\fr{1}{d(x)^2}|b_1(y)|(d(x) + |x-y|)e^{-\fr{|x-y|^2}{Kt}}dydt\\
\le Cd(x)^{-1}\int_{|x-y|\ge d(x)} |b_1(y)||x-y|^{-d+1}(1 + \fr{|x-y|}{d(x)})dy\\
\le  Cd(x)^{-\fr{d}{p}}\|b_1\|_{L^p}.
\ea
\la{nauoutoutb}
\ee
We used for the last integral $d>2$. In the case $d=2$ we take advantage of the fact that $H_D$ is smooth for large time to bound
\be
\ba
d(x)^{-2}\int_0^T t^{-\fr{3}{2}}dt\int_{|x-y|\ge d(x)} |x-y|e^{-\fr{|x-y|^2}{Kt}}|b_1(y)|dy\\
\le CTd(x)^{-2}\int_{|x-y|\ge d(x)} |x-y|^{-2}|b_1(y)|dy\\
 \le Cd(x)^{-\fr{2}{p}}\|b_1\|_{L^p}.
\ea
\la{naoutoutb2}
\ee
We have thus
\be
|\na u^{out}(x)|\le C_{\epsilon}d(x)^{-\fr{d}{p}}\left (1 + \left(\fr {d(x)}{\rho}\right)^{d+1}\right)\|b_1\|_{L^{p}(\Omega)}
\la{nauoutb}
\ee
with $C$ depending on $p$, $\epsilon$ and $\Omega$ only. 

We split
\be
\delta_h u = \delta_h u^{in} + \delta_h u^{out}
\la{split}
\ee
with
\be
\delta_h u(x)^{in} = c\int_0^{\rho^2}t^{-\fr{1}{2}}dt\int_{\Omega}\delta_h^x\na_
x^{\perp}H_D(x,y,t)\theta(y)dy
\la{duhin}
\ee
with $\rho$ satisfying $\rho\le \epsilon d(x)$. From the bound (\ref{nauoutb}) we have
\be
|\delta _h u^{out}(x)| \le C_{\epsilon}|h|d(x)^{-\fr{d}{p}}\left (1 + \left(\fr {d(x)}{\rho}\right)^{d+1}\right)\|b_1\|_{L^{p}(\Omega)}
\la{deltahuoutbg}
\ee
with a constant depending only on $\Omega$, $\epsilon$  and $p$. Note that if 
\be
\rho(x) = \epsilon d(x)
\la{rhoeq}
\ee
then  the estimate becomes
\be 
|\delta _h u^{out}(x)| \le C_{\epsilon}|h|d(x)^{-\fr{d}{p}}\|b_1\|_{L^{p}(\Omega)}
\la{deltahuoutb}
\ee
with a constant depending only on $\Omega$, $\epsilon$  and $p$.

We take a standard cutoff $\phi$ with scale $\ell$ and its companion $\chi$ and write
\be
\phi(x)\delta_h u^{in}(x) = u_h(x) + v_h(x)
\la{deltahu}
\ee
with 
\be 
u_h(x) = c\int_0^{\rho^2}t^{-\fr{1}{2}}dt\int_{\Omega}\na_x^{\perp}H(x,y,t)(\chi(y)\delta_h\theta(y)- \chi(x)\delta_h\theta(x))dy 
\la{uh}
\ee
and where
\be
v_h(x) = e_1(x) + e_2(x) + e_3(x) + e_5(x) + \phi(x)\delta_h\theta(x)e_4(x)
\la{vh}
\ee
with
\be
e_1(x) = c\int_0^{\rho^2}t^{-\fr{1}{2}}dt\int_{\Omega}\na_x^{\perp}(H_D(x+h,y,t)-H_D(x,y,t))\phi(x)(1-\chi(y))\theta(y)dy,
\la{e1}
\ee
\be
e_2(x) =  c\int_0^{\rho^2}t^{-\fr{1}{2}}dt\int_{\Omega}\na_x^{\perp}(H_D(x+h,y,t)-H_D(x,y-h,t))\phi(x)\chi(y)\theta(y)dy,
\la{e2}
\ee
\be
e_3(x) =c\int_0^{\rho^2}t^{-\fr{1}{2}}dt\int_{\Omega}\na_x^{\perp}H_D(x,y,t)(\chi(y+h)-\chi(y))\phi(x)\theta(y+h)dy,
\la{e3}
\ee
and
\be
e_4(x) =c\int_0^{\rho^2}t^{-\fr{1}{2}}dt\int_{\Omega}\na_x^{\perp}H_D(x,y,t)dy.
\la{e4}
\ee
We used here the facts that $\phi = \chi\phi$ and that
$(\chi\theta)(\cdot)$ and  $(\chi\theta)(\cdot + h)$
are compactly supported in $\Omega$ and hence 
\[
\int_{\Omega}\na_x^{\perp}H_D(x,y-h,t)\phi(x)\chi(y)\theta(y)dy = \int_{\Omega}\na_x^{\perp}H_D(x,y,t)\phi(x)\chi(y+h)\theta(y+h)dy.
\]

From (\ref{naxnaxb}) we have
\[
|e_1(x)|\le C\rho^2(x)|h|\int_0^1d\lambda \int_{A}\fr{1}{|x +\lambda h -y |^{d+3}}(d(x) + |x-y|)|b_1(y)|dy
\]
where $A = \{y\in\Omega \left |\right.\; |x_0-y|\ge \fr{5\ell}{8}\}$ is the support of $1-\chi(y)$. Because $x$ belongs to the support of $\phi$, it follows that
\be
|e_1(x)| \le C |h|d(x)^{-\fr{d}{p}}\|b_1\|_{L^p(\Omega)}.
\la{e1b}
\ee
We bound $e_3$ using $|\na\chi|\le C\ell^{-1}$, Lemma \ref{intpk} with $m= d$ in conjunction with (\ref{grbx}):
\[
|e_3(x)|\le C{|h|}\int_{\Omega}|\na\chi(y)|\fr{1}{|x-y|^d}|b_1(y)|dy,
\]
and consequently, because $x$ belongs to the support of $\phi$,
\be
|e_3(x)| \le C |h|d(x)^{-\fr{d}{p}}\|b_1\|_{L^p(\Omega)}.
\la{e3b}
\ee
Regarding $e_4$, in view of
\be
\int_{\Omega}\na_y^{\perp}H_D(x,y,t)dy = 0
\la{intnah}
\ee
we have
\[
e_4(x) = c\int_0^{\rho^2}t^{-\fr{1}{2}}\int_{\Omega}\left(\na_x^{\perp} + \na_y^{\perp}\right)H_D(x,y,t)dy.
\]
From (\ref{cancel1}) and Lemma \ref{intpk} with 
$m=j=0$, choosing  $\epsilon = \epsilon(\delta)$ small enough in $\rho= \epsilon d(x)$, we obtain
\be
|e_4(x)| \le \delta.
\la{e4b}
\ee
In order to estimate $e_2$ we write
\be
H_D(x+h,y,t)- H_D(x,y-h,t) = h\cdot\int_0^1 (\na_x+\na_y)H_D(x+\lambda h, y + (\lambda-1)h,t)d\lambda
\la{transh}
\ee
and use (\ref{cancel21})  to obtain
\[
\ba
|e_2(x)|\\
\le |h|\int_0^1d\lambda \int_0^{\rho^2}t^{-\fr{1}{2}}dt\int_{\Omega}|\na_x^{\perp}(\na_x +\na_y)H_D(x+\lambda h, y+(\lambda-1)h)||\chi(y)|(d(x) + |x-y|)|b_1(y)|dy\\
\le C|h|\|b_1\|_{L^{p}}d(x)^{1+ d(1-\fr{1}{p})}\int_0^{\rho^2}t^{-\fr{d+3}{2}}e^{-\fr{d(x)^2}{4Kt}}dt,
\ea
\]
and therefore
\be
|e_2(x)| \le C |h|d(x)^{-\fr{d}{p}}\|b_1\|_{L^p(\Omega)}.
\la{e2b}
\ee
Finally, by a Schwartz inequality,
\be
|u_h(x)| \le C\sqrt{\rho D(\chi\delta_h\theta)}.
\la{uhb}
\ee
This concludes the proof of the proposition.

We give also a bound for the normal component of the velocity at the boundary.
\beg{prop}\la{unb} Let $T$ be a $C^1$  divergence-free vector field tangent to the boundary of $\Omega$. Let $N = -T^{\perp}$. There exist constants $\ell_0>0$ and $C$ depending on the domain $\Omega$ such that
\be
|u(x)\cdot N(x)| \le C\left(d(x)^{1-\fr{d}{p}}\|b_1\|_{L^p} + d(x)^{\alpha}\|\theta\|_{C^{\alpha}(\Omega)}\right)\|T\|_{L^{\infty}} + Cd(x)^{2-\fr{d}{p}}\|b_1\|_{L^p}\|\na T\|_{L^{\infty}}
\la{unbound}
\ee
holds for $d(x)\le \ell_0$.
\end{prop}    
\noindent {\bf Proof.} In view of (\ref{nauoutb}), the fact that $T = N^{\perp}$ is tangent to the boundary, and $T\cdot\na_xH_D(x,y,t) =0$ for $x\in\pa\Omega$, we have
\be
|u^{out}(x)\cdot N(x)| \le C d(x)^{1-\fr{d}{p}}\|b_1\|_{L^p}
\la{uoutnb}
\ee
where 
\be 
u^{out}(x) = \int_{cd^2(x)}^{\infty} t^{-\fr{1}{2}}dt\int_{\Omega}\na_x^{\perp}H_D(x,y,t)\theta(y)dy.
\la{uout}
\ee
We consider now $u^{in}$  and  we write
\be
\ba
u^{in}(x)\cdot N(x) = -\int_0^{cd^2(x)}t^{-\fr{1}{2}}dt \int_{\Omega}T(x)\cdot\na_x H_D(x,y,t)\theta(y)dy\\
= \int_0^{cd^2(x)}t^{-\fr{1}{2}}dt\int_{\Omega} (T(y)-T(x))\cdot\na_x H_D(x,y,t)\theta(y)dy \\- \int_0^{cd^2(x)}t^{-\fr{1}{2}}dt\int_{|x-y|\ge cd(x)} T(y)\cdot(\na_x +\na_y) H_D(x,y,t)\theta(y)dy\\
-\int_0^{cd^2(x)}t^{-\fr{1}{2}}dt\int_{|x-y|\le cd(x)} T(y)\cdot(\na_x + \na_y) H_D(x,y,t)\theta(y)dy \\ + \int_0^{cd^2(x)}t^{-\fr{1}{2}}dt \int_{\Omega}T(y)\cdot\na_y H_D(x,y,t)(\theta(y)-\theta(x))dy\\ = U_1+U_2 +U_3 + U_4.
\ea
\la{unin15}
\ee
Because $|T(x)-T(y)|\le C|x-y|$ we have that
\[
|U_1| \le \|\na T\|_{L^{\infty}}\int_0^{cd^2(x)}t^{-\fr{1}{2}}dt\int_{\Omega}t^{-\fr{d}{2}}\fr{|x-y|}{t^{\fr{1}{2}}}e^{-\fr{|x-y|^2}{Kt}}(d(x) + |x-y|)|b_1(y)|dy
\]
and therefore
\be
|U_1| \le C d(x)^{2-\fr{d}{p}}\|b_1\|_{L^{p}}\|\na T\|_{L^{\infty}}
\la{U1b}
\ee
holds. For $U_2$ we use the bounds (\ref{grbx}), (\ref{grby}), $\theta = w_1 b$ (\ref{triangle}) and  Lemma \ref{intpk} to obtain
\be
|U_2| \le C\int_{|x-y|\ge d(x)}|x-y|^{-d}\left(d(x) + |x-y|\right)|b_1(y)|dy \le
Cd(x)^{1-\fr{d}{p}}\|b_1\|_{L^p}
\la{U2b}
\ee
with the understanding that if $p=\infty$ then $d(x)^{1-\fr{d}{p}}$ is replaced by $d(x)\log\left(\fr{\text{diam}\, \Omega}{d(x)}\right)$.  
For $U_3$ we use the bound (\ref{cancel12}) and Lemma \ref{intpk} 
to write
\be
d(x)\int_0^{cd^2(x)}t^{-\fr{1}{2}}\int_{|x-y|\le d(x)}|(\na_x +\na_y)H_D(x,y,t)||b_1(y)|dy
\le Cd(x)^{1-\fr{d}{p}}\|b_1\|_{L^p}.
\la{U3b}
\ee
Finally, for $U_4$ we use (\ref{grby}) and the fact that
\be
\int_0^{cd^2(x)}t^{-1-\fr{d}{2}}\int_{\Rr^d}|x-y|^{\alpha}e^{-\fr{|x-y|^2}{t}}dy\le Cd(x)^{\alpha}
\la{elem}
\ee
to obtain
\be
|U_4| \le Cd(x)^{\alpha}\|\theta\|_{C^{\alpha}}.
\la{U4b}
\ee
This concludes the proof of the lemma.
\section{Commutators}\la{commutator}
We consider the finite difference
\be
(\delta_h\l\theta)(x) = (\l\theta)(x+h)-(\l\theta)(x)
\la{dht}
\ee
with  $|h|\le \fr{d(x)}{32}$. We use a standard cutoff with scale $\ell$, $\phi$ and its companion $\chi$. 

\beg{prop}\la{commuh} We consider the commutator
\be
C_h(\theta) = \phi(x)(\delta_h\l\theta)(x) - \phi(x)\l(\chi\delta_h\theta)(x).
\la{ch}
\ee
There exists a constant $\Gamma_0$ such that the  commutator $C_h(\theta)$ obeys

\be
\left |C_h(\theta)(x)\right| \le \Gamma_0\fr{|h|}{d(x)}\|b_1\|_{L^{p}(\Omega)}
d(x)^{-\fr{d}{p}}
\la{commhb}
\ee
for $|h|\le\fr{\ell}{16}$, $\theta\in H_0^1(\Omega)\cap L^{\infty}(\Omega)$ and $b_1 = \fr{\theta}{w_1}\in L^p(\Omega)$ with $p>{d}$. The constant is bounded as $p\to \infty$ and if $p=\infty$ the estimate is
\be
\left |C_h(\theta)(x)\right| \le \Gamma_0\fr{|h|}{d(x)}\|b_1\|_{L^{\infty}(\Omega)}.
\la{commhbb}
\ee

\end{prop}
\noindent{\bf Proof.} 
We compute the commutator as follows   
\be
\ba
(\phi \delta_h\l\theta)(x) - \phi(\l \chi\delta_h\theta)(x)\\ = c\int_0^{\infty}t^{-\fr{3}{2}}dt\int_{\Omega}(H_D(x,y,t)-H_D(x+h,y,t))\phi(x)(1-\chi(y))\theta(y)dy\\
-c\int_0^{\infty}t^{-\fr{3}{2}}dt\int_{\Omega}(H_D(x+h,y,t)-H_D(x,y-h,t))\phi(x)\chi(y)\theta(y)dy\\
-c\int_0^{\infty}t^{-\fr{3}{2}}dt\int_{\Omega}H_D(x,y,t)\phi(x)(\delta_h\chi)(y)\theta(y+h)dy \\
 =  E_1(x) + E_2(x) + E_3(x). 
\ea
\la{compucom}
\ee

We use (\ref{compucom}).  
We  observe by triangle inequlaity $d(y) \le d(x) + |x-y|$  and thus
\be
|\theta(y)| \le C|b_1(y)|(d(x) + |x-y|)
\la{thetaby}
\ee
holds for any $x,y\in \Omega$.
For $E_1(x)$ we use  the inequalities (\ref{hb}), (\ref{grbx}),  and Lemma \ref{intpk} with $m=d+2$ when $t\le d(x)^2$,  and $m=d+1$ when $t\ge d(x)^2$, together with $d(x)^{-1}H_D \le C|x-y|^{-1}t^{-\fr{d}{2}}e^{-\fr{|x-y|^2}{Kt}}$. Substituting  (\ref{thetaby}) for $\theta$, we deduce
\[
|E_1(x)|\le C|h|\int_{\Omega}[|x-y|^{-(d+2)}(d(x) +|x-y|)|b_1(y)|\phi(x)|1-\chi(y)|dy
\]
and then, using a H\"{o}lder inequality we obtain
\be
E_1(x) \le C \fr{|h|}{d(x)}\|b_1\|_{L^p}d(x)^{-\fr{d}{p}}.
\la{E1b}
\ee
For $E_2$ we use (\ref{transh}) like in the proof of the estimate (\ref{e2b})
and (\ref{cancel12}),  together with Lemma \ref{intpk} with $m=d+2$, $j=0$, 
and bounds 
\be
\ba
\int_{0}^{\infty}t^{-\fr{3}{2}}t^{-\fr{d+1}{2}}e^{-\fr{d(x0^2}{Kt}}dt\int_{|x-y|\le d(x)}|b_1(y)|(d(x) + |x-y|) dy\\
\le C(d(x))^{-2-d}d(x)\int_{|x-y|\le d(x)}|b_1(y)|dy\\
 \le Cd(x)^{-1-\fr{d}{p}}\|b_1\|_{L^p}
\ea
\la{naxnaylongtimein}
\ee
and
\be
\ba
\int_0^\infty t^{-\fr{3}{2}}\int_{|x-y|\ge d(x)}(t^{-\fr{1}{2}} + \fr{1}{d(x)} + \fr{1}{d(y)})H_D(x,y,t) |b_1(y)|(d(x) + |x-y|)dydt\\\le
C\int_0^\infty t^{-\fr{3}{2}}\int_{|x-y|\ge d(x)}(t^{-\fr{1}{2}} + \fr{1}{|x-y|})t^{-\fr{d}{2}}e^{-\fr{|x-y|^2}{Kt}}|b_1(y)|(d(x) + |x-y|)dy\\
\le C\int_{|x-y|\ge d(x)}\fr{1}{|x-y|^{d+2}}|b_1(y)|(d(x) + |x-y|)dy\le Cd(x)^{-1-\fr{d}{p}}\|b_1\|_{L^p}
\ea
\la{naxnayalltimenout}
\ee
to obtain
\be
|E_2(x)|\le  C\fr{|h|}{d(x)}\|b_1\|_{L^{p}}{d(x)}^{-\fr{d}{p}}.
\la{E2b}
\ee

For $E_3$ we have
\[
|E_3(x)|\le |h|\int_0^{\infty}t^{-\fr{3}{2}}dt\int_{\Omega}H_D(x,y,t)\phi(x)|\na \chi(y)|(d(x) + |x-y|)b_1(y)dy
\]
and from Lemma \ref{intpk} with $m=d, d+1$, $j=0$ we obtain
\[
|E_3(x)|\le  C\fr{|h|}{d(x)}\|b_1\|_{L^{p}}d(x)^{-\fr{d}{p}}.
\]

This concludes the proof.



\section{SQG: H\"{o}lder bounds}\la{holder}
We consider the equation (\ref{sqg}) with $u$ given by (\ref{u}) and with initial data $\theta_0\in H_0^1(\Omega)\cap L^{\infty}(\Omega)$.  We note  we have
\be
\|\theta(t)\|_{L^{\infty}}\le \|\theta_0\|_{L^{\infty}}.
\la{linf}
\ee
We prove the following result. 
\beg{thm}\la{glh} Let $\theta(x,t)$ be a solution of (\ref{sqg}) in the bounded domain with smooth boundary $\Omega$, obeying (\ref{gradintb}) on a time interval $[0,T]$. Assume that 
\be
\sup_{0\le t\le T}\|b_1(t)\|_{L^p(\Omega)} \le B
\la{b1B}
\ee
holds with $p>d$. Then for $0<\alpha<1-\fr{d}{p}$ there exists a constant $K$,  depending only on the domain $\Omega$ and p, such that 
\be
\sup_{0\le t\le T}\sup_{x\in\Omega}\sup_{|h|\le\fr{d(x)}{32}}\fr{|\delta_h \theta(x,t)|}{|h|^{\alpha}} \le 2\|\theta_0\|_{C^{\alpha}} + KB(M+1) 
\la{hldrbnd}
\ee
holds, where $M$ is the a priori bound in (\ref{gradint}). Moreover, the velocity $u$ is bounded $u\in L^{\infty}(0,T; L^{\infty}(\Omega))$, obeying (\ref{ucondb}) and the normal component of the velocity vanishes near the boundary of order $d(x)^{\alpha}$, obeying (\ref{unbound}).

\end{thm}

\noindent{\bf{Proof.}} 
We take
\be
|h|\le\fr{\ell}{16}.
\la{hdeltal}
\ee
We take $x$ with $d(x)\le 2\ell_0$. From the SQG equation we obtain the equation
\be
\fr{1}{2}\left(\pa_t + u\cdot\na\right)|\delta_h\theta|^2 + (\delta_h\theta)\delta_h\l\theta = -(\delta_h\theta) \delta_h u\cdot \na \theta(x+h).
\la{sqgsq}
\ee 
We use a standard cutoff $\phi$ with scale $\ell$, and companion $\chi$. 
We multiply by $\phi^2$ and obtain
\be
\fr{\phi^2}{2}(\pa_t + u\cdot\na )(|\delta_h\theta|^2) + \phi^2 \delta_h\theta\l(\chi \delta_h\theta) = - (\phi \delta_h\theta) C_h(\theta) - (\phi \delta_h\theta)\phi \delta_h u\cdot \na\theta(x+h) 
\la{le}
\ee
where $C_h(\theta)$ is the commutator given above in (\ref{ch}). 

Multiplying by $|h|^{-2\alpha}$ where $\alpha>0$ is smaller than $1-\fr{d}{p}$,  we obtain
\be
\fr{\phi^2}{2}\left(\pa_t + u\cdot\na\right)(f^2) + \phi^2 f\l(\chi f)  = 
-|h|^{-\alpha}\phi f C_h(\theta)   - |h|^{-\alpha}\phi f \phi\delta_hu \cdot \na\theta(x+h)
\la{bu}
\ee
where 
\be
f(x,t; h) = f = |h|^{-\alpha}\delta_h\theta (x,t).
\la{f}
\ee
The first term in the right hand side of (\ref{bu}) is bounded using the commutator estimate (\ref{commhb}). 
\be
\ba
 \phi |f||h|^{-\alpha} |C_h(\theta)|\le 
\phi|f||h|^{-\alpha}\Gamma_0\fr{|h|}{d(x)}Bd(x)^{-\fr{d}{p}}\\
\le (C\Gamma_0 B |h|^{1-\alpha-\fr{d}{p}}) \fr{1}{d(x)} \phi |f|
\\
\le \fr{1}{d(x)}\left [C\Gamma_0 B \ell^{1-\alpha-\fr{d}{p}}\right]\phi |f|.
\ea
\la{commbo}
\ee
In view of (\ref{gradintb}) we have that
\be
|\na \theta(x+h)|\le M\fr{1}{d(x)}.
\la{deltahm}
\ee
The second term in the right hand side of (\ref{bu}) is estimated using (\ref{deltauhb}) with $\delta =\delta(\epsilon) \le \fr{\gamma_1}{8M}$  with $\epsilon$ sufficiently small (depending on $\gamma_1$ and $M$ but not on $B$). We obtain 
\be
\ba
\phi |f||h|^{-\alpha}|\phi\delta_h u| |\na\theta(x+h)|\\
\le
Md(x)^{-1}\phi |f||h|^{-\alpha}\left[\sqrt{\epsilon d(x)D(\chi \delta_h\theta)} + C_{\epsilon}|h|d(x)^{-\fr{d}{p}}B +\fr{\gamma_1}{8M}\phi|\delta_h\theta(x)|\right]\\
\le \fr{1}{2}D(\chi f) + \fr{\gamma_1}{4d(x)}\phi^2 |f|^2 +
C_{\epsilon}BM|h|^{1-\alpha}d(x)^{-1-\fr{d}{p}}\phi|f|.
\ea
\la{third}
\ee
where $D(g)$ is given in (\ref{dfdxb}) and where we also used $M^2\epsilon\le \fr{\gamma_1}{8}$. Therefore, if we have
\be
0<\alpha< 1- \fr{d}{p}
\la{alphapbound}
\ee
we obtain from (\ref{bu}), (\ref{commbo}), (\ref{third}) that
\be
\fr{\phi^2}{2}\left(\pa_t + u\cdot\na\right)(f^2) + \phi^2 f\l(\chi f) \le 
\fr{1}{2}D(\chi f) + \fr{\gamma_1}{4d(x)}\phi^2 |f|^2 + \fr{1}{d(x)}\left [K_1B(M+1) \ell^{1-\alpha-\fr{d}{p}}\right]\phi |f|
\la{fin}
\ee
holds for $|h|\le\fr{\ell}{16}$. Note that $K_1$ does not depend on $\ell$ nor on $h$ and that, in view of (\ref{alphapbound}) we may take $|h|$ and $\ell >0 $ as small as we wish.

The rest of the argument is by contradiction. We fix $T>0$ and take $0<\ell<\ell_0$. We consider the compact region
\be
A_{\ell} = \{x\in \Omega\left |\right.\; \ell\le d(x)\le 2\ell\}.
\la{al}
\ee
We assume by contradiction that there exists $x_1\in A_{\ell}$, $t_0\in [0,T)$ and $h_0$ with $|h_0|\le \fr{d(x_1)}{32}$ such that
\be
|h_0|^{-\alpha}|{\delta_{h_0}\theta(x_1, t_0)}| \ge 2 \|\theta_0\|_{C^{\alpha}} + 
KB(M+1)\ell_0^{1-\alpha -\fr{d}{p}}.
\la{contra}
\ee
with $K = \fr{5K_1}{\gamma_1}$ where $K_1$ appears in (\ref{fin}). We may assume without loss of generality that $t_0$ is the infimum of such $t$ that (\ref{contra}) holds for some $x_1\in A_{\ell}$. The prefactor $2$ in front of $\|\theta_0\|_{C^{\alpha}}$ was put there for convenience, in order to make sure that $t_0>0$ (it could have been any number larger than one). We fix $h_0$ and $t_0$ and take $x_0\in A_{\ell}$ to be a point where the maximum of the function $f^2(x, t_0; h_0)$ is achieved in the region $A_{\ell}$. We know that $\theta$ is interior Lipschitz, so $|h_0|>0$ and $f^2$ is Lipschitz continuous there. Therefore (\ref{contra}) holds with $x_0$ replacing $x_1$. We take a standard cutoff with scale $\ell$, center $x_0$ and companion $\chi$.
We use the inequality (\ref{cor}):
\be
\chi f\l (\chi f) - \fr{1}{2}\l(\chi^2f^2) = D(\chi f) \ge  \gamma_1 (d(x))^{-1} \chi^2|f|^2,
\la{d1}
\ee
valid pointwise.   We also use the fact that $\l (\chi^2f^2)(x_0)>0$ because  $\chi^2f^2$ is maximized at $x_0$.  In fact, more is true,
\[
\l(\chi^2 f^2)(x_0) \ge \chi^2(x_0)f^2(x_0)\l 1.
\]
Indeed, for any function $g$ in the domain of $\l$ which achieves its maximum at $x_0\in \Omega$ we have
\[
(\l g)(x_0) = c\int_0^{\infty}t^{-\fr{3}{2}}\left(g(x_0)-\int_{\Omega}H_D(x,y,t)g(y)dy\right)\ge g(x_0)\l 1.
\]
Using $\phi(x_0) = \chi(x_0) =1$ we have from (\ref{fin}), (\ref{contra}) and (\ref{d1}) and the fact that $u\cdot\na f^2 = 0$ at an interior local maximum, 
\be
\ba
\fr{1}{2}\pa_t f^2  \le  - \chi f\l\chi f +\fr{1}{2}D(\chi f) + \fr{\gamma_1}{4d(x_0)}\phi^2 f^2 + \fr{1}{d(x_0)}K_1B(M+1)\ell^{1-\alpha-\fr{d}{p}}\phi|f| \\
\le -\fr{1}{2}\l(\chi^2 f^2) - \fr{1}{2}D(\chi f) + \fr{\gamma_1}{4d(x_0)}\phi^2 f^2 +  \fr{1}{d(x_0)}K_1B(M+1)\ell^{1-\alpha-\fr{d}{p}}\phi|f|\\
\le -\fr{\gamma_1}{4d(x_0)}\phi^2 f^2 +  \fr{1}{d(x_0)}K_1B(M+1)\ell^{1-\alpha-\fr{d}{p}}\phi|f|\\
\le -\fr{\gamma_1}{4d(x_0)}B(M+1)\ell_0^{1-\alpha-\fr{d}{p}}\phi|f|(K-\fr{4K_1}{\gamma_1})\\
\le -\fr{1}{4d(x_0)}B(M+1)\ell_0^{1-\alpha-\fr{d}{p}}K_1<0,
\ea 
\la{fina}
\ee
which is a contradiction.
Thus (\ref{hldrbnd}) holds, and the proof of the uniform bound on the H\"{o}lder norm is concluded. The fact that $u$ obeys (\ref{ucondb}) follows from Lemma \ref{uconditional} and the vansihing of the normal component of velocity follows from Proposition \ref{unb}, in view of the bound (\ref{hldrbnd}).


{\bf{Acknowledgment.}} The work of PC was partially supported by  NSF grant DMS-1209394

\end{document}